\documentclass[11pt]{amsart}
\usepackage{amssymb}
\usepackage{verbatim} 
\newtheorem{thm}[equation]{Theorem}
\newtheorem{pro}[equation]{Proposition}
\newtheorem{cor}[equation]{Corollary}
\newtheorem{lem}[equation]{Lemma}
\newtheorem{exa}[equation]{Example}

\newtheorem{deft}[equation]{Definition}

\numberwithin{equation}{section}

\def\a1s{a_1,\cdots, a_s}
\def\a{\alpha}

\def\aa{\mathcal A}

\def\andd{\quad\hbox{and}\quad}

\def\fb{\frak{b}}
\def\b{\beta}

\def\bb{\mathcal{B}}

\def\bl4{B_{\ell\geq4}}
\def\cc{{\mathcal C}}

\def\d{\delta}
\def\D{\Delta}

\def\dd{\mathcal D}

\def\bbbf{\mathbb{F}}

\def\gg{{\mathcal G}}

\def\fg{\mathfrak{g}}

\def\hh{{\mathcal H}}

\def\Lam{\Lambda}
\def\LL{\mathcal{L}}

\def\LLc{\mathcal{L}_c}
\def\LLcc{\mathcal{L}_{cc}}

\def\m{\mathcal{M}}
\def\bbbn{\mathbb{N}}

\def\1k{\frac{1}{k}}
\def\op{\oplus}
\def\ot{\otimes}

\def\qed{\hfill$\Box$}
\def\sub{\subseteq}

\def\pf{\noindent{\bf Proof. }}

\def\ss{\mathcal{S}}

\def\T{{\mathcal T}}

\def\v{{\mathcal V}}

\def\w{{\mathcal W}}

\def\bbbz{{\mathbb Z}}
\def\ss{\mathcal{S}}
 
\def\1il{1\leq i\leq\ell}

\begin{document}
\title{Locally finite admissible  simple Lie algebras }

\author{ Malihe Yousofzadeh}
\address
{Department of Mathematics\\ University of Isfahan\\Isfahan, Iran,
P.O.Box 81745-163} \email{ma.yousofzadeh@sci.ui.ac.ir.}
\keywords{\em Locally finite Lie algebras,  root graded
Lie algebras.\\{2010 Mathematics Subject Classification(s):17B05, 17B20,
17B60,  17B70, 17B67}
}

\begin{abstract}
We introduce a class of Lie algebras called admissible Lie algebras. We show that a locally finite admissible simple Lie algebra contains a nonzero maximal toral subalgebra and the corresponding root system is an irreducible locally finite root system.
\end{abstract}
\maketitle
\section{introduction}
In 1976, G.B. Seligman \cite{Se} showed that finite dimensional simple Lie algebras containing a nonzero  maximal toral subalgebra have a decomposition as $$(*)\hspace{1cm}\LL=(\gg\ot \aa)\op(\ss\ot \bb)\op(\v\ot \cc)\op\dd$$ in which $\gg$ is a finite dimensional split  simple Lie algebra and
$\v,$ $\ss$ are specific irreducible $\gg$-modules. The Lie algebra structure on $\LL$ induces an  algebra structure on   the vector space $\fb:=\aa\op\bb\op\cc$ and $\dd$ is a Lie subalgebra  of $\LL$ isomorphic to a specific subalgebra of inner derivations of $\fb.$ Following G.B. Seligman, S. Berman and R. Moody  \cite{BM} introduced the notion of a Lie algebra graded by an irreducible  reduced finite root system $R$ which  is  a Lie algebra containing a finite dimensional split  simple Lie subalgebra $\gg$ of type $R$  with a Cartan subalgebra $\hh$  such that $\LL$ has a weight space decomposition with respect to $\hh$ and as an algebra, $\LL$ is  generated by weight spaces corresponding to "nonzero" weights.
Next Allison, Benkart and Gao generalized \cite{BM}'s definition    to non reduced case by letting $\gg$ be a finite dimensional split simple Lie algebra of type $B,C$ or $D$ if $R$ is of type $BC.$
In  \cite{BM}, \cite{BZ} and \cite{ABG}, the authors stated recognition theorems for root graded Lie algebras to classify such  Lie algebras up to centrally  isogeny, indeed any root graded Lie algebra $\LL$ has a decomposition as $(*)$ with a prescribed  algebra structure on the so called coordinate algebra $\fb:=\aa\op\bb\op\cc.$  In 1999, N. Stumme  studied  locally finite  Lie algebras containing a splitting Cartan subalgebra and called such algebras, locally finite split  Lie algebras. She  showed that, a reduced locally finite root system,  the root system of  a locally finite split semisimple Lie algebra, is  the direct union of finite root subsystems  of semisimple types. She also proved that a locally finite split simple Lie algebra is the direct union of finite dimensional simple subalgebras.  Locally finite root systems  appear in the theory of Kac-Moody Lie algebras as well, more precisely, countable irreducible  reduced locally finite root systems are the root systems of infinite rank affine  Lie algebras
\cite[\S 7.11]{K}.  In 2003, O. Loos and E. Neher \cite{LN} introduced    locally finite root systems axiomatically and gave a complete description   of these root systems.  Next Neher \cite{N} generalized the definition of a root graded Lie algebra to Lie algebras graded by a  locally finite root system. In Neher's sense, if $R$ is a locally finite root system, $S$ is a subsystem of $R$ and $\Lam$ is an abelian group,  an $(R,S,\Lam)$-graded Lie algebra is a compatible  $\hbox{span}_\bbbz(R)$-graded and $\Lam$-graded Lie algebra whose support with respect to $\hbox{span}_\bbbz(R)$-grading is contained in $R,$ for  every $0\neq \a\in S,$ the homogeneous space $\LL_\a^0$
contains a so called invertible element and $\LL_0=\sum_{0\neq\a\in R}[\LL_\a,\LL_{-\a}].$ In this work, we first study  the direct union $\cup_{n\in\bbbn}\gg_n$ of finite dimensional simple Lie algebras $\gg_n,$ $n\in\bbbn,$ containing a nonzero maximal toral  subalgebra, these Lie algebras are  $(R,R_{sdiv},0)$-graded Lie algebras (see Definition \ref{def-root}) in Neher's sense for a locally finite root system $R.$ Here we study a class of Lie algebras, called  admissible Lie algebras. An admissible Lie algebra  $\LL$ has  a nonzero toral subalgebra $\hh$ contained in the subalgebra of $\LL$ generated by the weight spaces   such that for any nonzero weight vector $x,$ there is a weight vector $y$ such that $[x,y]\in\hh$ and  $(x,[x,y],y)$ is an $\mathfrak{sl}_2$-triple. We show that if the Lie algebra $\LL$ is locally finite simple Lie algebra whose weight spaces are finite dimensional, then $\hh$ is a maximal toral subalgebra, the root system of $\LL$ with respect to $\hh$ is a locally finite root system and $\LL$ is the direct union of finite dimensional simple subalgebras.
\section{admissible Lie algebras}
Throughout   this work $\bbbn$ denotes the set of all nonnegative
integers and $\bbbf$ is a field of characteristic zero. Unless
otherwise mentioned, all vector spaces are considered over
$\bbbf.$ In the present paper, we denote the dual space of a
vector space $V$ by $V^\star$ and by $GL(V),$ we mean the group of
automorphisms of $V.$
For a matrix $A,$ $tr(a)$ denotes the trace of $A.$ Also for a Lie algebra $\LL,$   we mean by $Z(\LL),$ the center of
$\LL$ and if $\LL$ is finite dimensional, we denote the Killing
form of $\LL$ by $\kappa.$ We also make a convention that for
elements  $x_1,\ldots,x_m$ of a Lie algebra, by an expression of
the form $[x_1,\ldots,x_m],$ we always mean $[x_1,\ldots[x_{m-1},
x_m]\ldots].$
to be zero.
\begin{deft}\label{def-root}
{\rm Let $\v$ be a nontrivial vector space and $R$ be a subset of
$\v,$ $R$ is said to be a {\it locally finite root system in $\v$}
if

(i) $0\not\in R,$ $R$ is locally finite and spans $\v,$

(ii) for every $\a\in R,$ there exists $\check\a\in\v^\star$ such
that $\check\a(\a)=2$ and $s_\a(\b)\in  R$ for $\a,\b\in  R$ where
$s_\a:\v\longrightarrow\v$ maps $v\in\v$ to $v-\check\a(v)\a.$

(iii) $\check\a(\b)\in\bbbz,$ for $\a,\b\in R.$

Set  $R_{sdiv}:=R\setminus\{\a\in R\mid 2\a\in R\}$ and call it  the semi-divisible  subset of $R,$ the root system $R$ is called {\it reduced} if $R=R_{sdiv}.$

}
\end{deft}

Suppose that $R$  is a locally finite root   system. A nonempty
subset $S$ of $R$  is said to be  a {\it subsystem} of $R$ if
$s_\a(\b)\in S$ for $\a,\b\in S.$ Following \cite[\S 2.6]{LN}, we
say two roots $\a,\b$ are {\it connected} if there exist finitely
many roots $\a_1=\a,\a_2,\ldots,\a_n=\b$ such that
$\check\a_{i+1}(\a_{i})\neq 0,$ $1\leq i\leq n-1.$ Connectedness
is an equivalence relation on $R.$  The root system $R$ is the
disjoint union of its connected components. A nonempty subset $X$
of $R$ is called {\it irreducible,} if each two elements $x,y\in
X$ are connected and it is called {\it closed} if for $x,y\in X$
with $x+y\in R,$ we would have $x+y\in X.$ It is easy to see that
each connected component of the locally finite root system $R$ is
a closed subsystem of $R.$ Also using \cite[Cor. 3.15]{LN}, $R$
is a direct limit of its finite subsystems, and if $R$ is
irreducible, it is a direct limit of its irreducible finite
subsystems.
\begin{deft}
{\em Let $\hh$ be a Lie algebra. We say an $\hh-$module $\m$ has a
{\it weight space decomposition with respect to $\hh,$} if
\begin{equation*}
\m=\op_{\a\in \hh^\star}\m_\a \;\hbox{where}\; \m_\a:=\{x\in\m\mid
h\cdot x=\a(h)x;\;\;\forall\; h\in \hh\};\;\a\in\hh^\star.
\end{equation*}
The set $R:=\{\a\in\hh^\star\setminus\{0\}\mid \m_\a\neq \{0\}\}$
is called the {\it set of weights} of $\m$ with respect to $\hh$
and $\m_\a,$ $\a\in R,$ is called a {\it weight space}, also any
element of $\m_\a$ is called a {\it weight vector} of {\it weight}
$\a.$ If a Lie algebra $\LL$ has a weight space decomposition with
respect to a subalgebra of $\LL$ via  the adjoint representation,
the set of weights is called the {\it root system } and weight
spaces are called {\it root spaces.}}
\end{deft}

An easy verification proves the following lemma.
\begin{lem}
\label{rest3} Let $\hh$ be a Lie algebra and $\m$ be  an
$\hh$-module admitting a weight space decomposition
$\m=\op_{\a\in\hh^\star}\m_\a$ with respect to $\hh.$ Let $T$ be a
subalgebra of $\hh$ and take $\pi:\hh^\star\longrightarrow T^\star$ to be
defined by $\pi(\a)=\a\mid_{_T},$ the restriction of $\a$ to $T.$
For $\b\in T^\star,$ define $\m'_\b:=\{v\in\m\mid t\cdot
v=\b(t)v;\;\forall t\in T\},$ then $\m=\op_{\b\in T^\star}\m'_{\b}$
and for $\b\in T^\star,$ $\m'_\b=\displaystyle{\bigoplus_{\{\a\in
\hh^\star,\;\pi(\a)=\b\}}\m_\a.}$
\end{lem}
\begin{deft}\label{split}
{\em  Let $\LL$ be a Lie algebra.  A nontrivial subalgebra $\hh$
of $\LL$ is called a {\it  split toral subalgebra}  if $\LL$ as an
$\hh-$module, via the adjoint representation, has a weight space
decomposition with respect to $\hh.$ One can see  that a split
toral subalgebra of a Lie algebra is abelian. Throughout the
present paper  by a toral subalgebra, we always mean a split toral
subalgebra.}
\end{deft}

Now suppose that $\LL$ is a Lie
algebra containing a nonzero toral subalgebra $\hh$ with
corresponding root system $R.$ The Lie algebra $(\LL,\hh)$ (or
$\LL$ for simplicity) is called {\it admissible } if $\LL$
satisfies the following property:

\begin{equation}\label{con}
\parbox{4.2in}{$\hh\sub\sum_{\a\in R}[\LL_\a,\LL_{-\a}]$ and for  $\a\in R,$  $0\neq x\in\LL_\a,$ there is $y\in\LL_{-\a}$
such that $h:=[x,y]\in\hh$ and $(x,h,y)$ is an
$\mathfrak{sl}_2$-triple. }
\end{equation}

An element $h\in\hh$ is called  a {\it splitting element
corresponding to $\a\in R$} if there are  $x\in\LL_\a$ and
$y\in\LL_{-\a}$ such that $(x, h:=[x,y],y)$ is an
$\mathfrak{sl}_2$-triple.  A subset $\D$ of $R$ is called {\it
connected} with respect to a fix set $\{h_\a\mid\a\in R\}$ of
splitting elements if for any $\b,\gamma\in\D,$ there is a finite
sequence $\a_1,\ldots,\a_n$ of elements of $\D$ such that
$\a_1=\b,$ $\a_n=\gamma$ and $\a_{i+1}(h_{\a_i})\neq 0,$ $1\leq
i\leq n-1.$  A root $\a\in R$ is called {\it integrable} if there
are $e_\a\in \LL_\a,$ $f_\a\in\LL_{-\a}$ such that
$h_\a:=[e_\a,f_\a]\in \hh,$ $(e_\a,h_\a,f_\a)$ is an
$\mathfrak{sl}_2$-triple and that $ad_{e_\a}$ and $ad_{f_\a}$ act
locally nilpotently on $\LL.$ We denote by $R_{int},$  the set of
integrable roots of $\LL$ and note that if $\LL$ is a locally
finite admissible Lie algebra, then $R=R_{int}.$ A subset $\D$ of
$R$ is called {\it symmetric } if $\D=-\D$ and it is called {\it
closed } if $(\D+\D)\cap R\sub \D.$

\begin{exa}\label{exa1}{\rm
Let $\LL$ be  a finite dimensional semisimple Lie algebra
containing a maximal toral subalgebra  $\T.$ Take $\Phi$ to be the
root system of $\LL$ with respect to $\T.$ Using \cite[Lem. I.3
and Lem. I.5]{Se}, we get that $(\LL,\T)$ is an admissible Lie
algebra. Moreover   Corollary to Lemma I.4 of \cite{Se} shows
that for $\b\in\Phi,$ there is a unique $k_\b\in\T$ with
$$
\b(k_\b)=2 \hbox{ and $k_\b=[u,v]$ for some
$u\in\LL_\b,v\in\LL_{-\b}.$}
$$
Using \cite[Lem. I.5 and Lem. I.6]{Se}, we get that
\begin{equation}\label{span}
\T=\hbox{span}_\bbbf\{k_\b\mid\b\in \Phi\}\andd
\T^\star=\hbox{span}_\bbbf\{\b\mid\b\in \Phi\}. \end{equation}
Following the proof of Lemma I.5 of \cite{Se}, we get that if
$\b\in\Phi,$  then  $\kappa(k_\b,k_\b)\neq 0$  and for $t\in\T,$
$\b(t)=2\kappa(t,k_\b)/\kappa(k_\b,k_\b),$ now (\ref{span})
implies that  the linear transformation form $\T$ to $\T^\star$
mapping $t\mapsto\kappa(t,\cdot),$ $t\in\T,$ is onto and so is one
to one, this in particular implies that the Killing form
restricted to $\T$ is non-degenerate and that for $\b\in\Phi,$
$2k_\b/\kappa(k_\b,k_\b)$ is the unique element of $\T$
representing $\b$ through the Killing form. Therefore one
concludes that
\begin{equation}\label{lin}
\parbox{4in}{\begin{center}if $\a,\b_1,\ldots,\b_n\in \Phi$ are such that $\a\in\hbox{span}_\bbbf\{\b_1,\ldots,\b_n\},$ then
$k_\a\in\hbox{span}_\bbbf\{k_{\b_1},\ldots,k_{\b_n}\}.$
\end{center} }
\end{equation}}
\end{exa}
\medskip

{\it From now on we assume $(\LL,\hh)$ is an admissible Lie
algebra with nonempty root system $R.$}

Using $\frak{sl}_2$-module  theory, one can easily prove the following proposition.
\begin{pro}
\label{loc5} Let $\a\in R.$ Take $e\in\LL_{\a}$ and
$f\in\LL_{-\a}$ to be such that $(e,h:=[e,f],f)$ is an
$\mathfrak{sl}_2$-triple and set
$\fg:=\hbox{span}_\bbbf\{e,h,f\}.$
Suppose    $ad_e$ and $ad_f$  act  locally nilpotently on
$\LL,$ then the followings are satisfied:

(i) For $\b\in R,$ $\b(h)\in\bbbz$  and $\b-\b(h)\a\in R.$

(ii) If $\b\in R$ is such that  $\a+\b\in R,$ we have
$[e,\LL_\b]\neq\{0\},$ in particular, $[\LL_\a,\LL_\b]\neq \{0\}.$

(iii) For $\b\in R,$ $\{k\in\bbbz\mid \b+k\a\in R\cup\{0\}\}$  is
an interval.

(iv) For $\b\in R,$ if $\b(h)>0,$ then $\b-\a\in R\cup\{0\}$ and if $\b(h)<0,$ then $\b+\a\in R\cup\{0\}.$

(v) If $k\in\bbbf$ and $k\a\in R,$ then $k\in \bbbz/2.$

(vi)  $\{k\a\mid k\in\bbbf\}\cap
R_{int}\sub\{\pm\a,\pm(1/2)\a,\pm2\a\}.$

\end{pro}

Now use  Proposition \ref{loc5}$(i),(iv)$ together with the same argument as in \cite[\S 10.2]{LN} to prove the following Corollary:

\begin{cor}\label{cor0}
If $R=R_{int,}$ $R$ has the partial sum property, in the sense that if $n\in\bbbn\setminus\{0\}$ and $\b,\a_1,\ldots,\a_n\in R\cup\{0\}$ are such that $\a_1+\cdots+\a_n=\b,$  there is a permutation $\pi$ of $\{1,\ldots,n\}$ such that $\a_{\pi(1)}+\cdots+\a_{\pi(i)}\in R\cup\{0\}$ for all $i\in\{1,\ldots,n\}.$
\end{cor}

\begin{pro}\label{semi-sim1} Suppose that $R=R_{int}.$

(a) For  an ideal $I$ of $\LL,$ set  $ R_{_I}:=\{\a\in R\mid
I\cap\LL_\a\neq \{0\}\},$ then  $ R_{_I}$ is a symmetric closed
subset of $ R$ and  $I= (I\cap\LL_0)\op \sum_{\a\in R}(I\cap
\LL_\a)=(I\cap\LL_0)\op \sum_{\a\in  R_{_I}}\LL_\a.$ Moreover if
$I$ is  simple  as   a Lie algebra, $I\cap\LL_0=\sum_{\a\in
R_{_I}}[\LL_\a,\LL_{-\a}].$
\smallskip

(b) If $\LL$ is semisimple (in the sense that it is a direct sum
of simple ideals), then $\LL_0=\sum_{\a\in R}[\LL_\a,\LL_{-\a}].$
\smallskip

(c) Let  $I$ be an ideal of $\LL.$ Define $\D:=R\setminus  R_{_I}$
and set $J:=\sum_{\a\in \D}[\LL_\a,\LL_{-\a}]\op\sum_{\a\in
\D}\LL_\a.$ If $\LL_0=\sum_{\a\in R}[\LL_\a,\LL_{-\a}],$ then $J$
is an ideal of $\LL,$ $\LL=I+J$ and $I\cap J\sub Z(\LL).$
\end{pro}

\pf $(a)$ Since $I$ is an ideal of $\LL,$ we have
$I=(I\cap\LL_0)\op\sum_{\a\in R}(I\cap\LL_\a)$ by \cite[Pro. 2.1.1]{MP}. Now
if $\a\in
 R_{_I}$ and $0\neq x\in I\cap\LL_\a,$ then by (\ref{con}), there
exists $y\in \LL_{-\a}$ such that $(x,h:=[x,y],y)$ is an
$\mathfrak{sl}_2$-triple. Now as $x\in I,$ one gets that $h,y\in
I$  and for all $z\in\LL_\a,$ $2z=\a(h)z=[h,z]\in I.$ So we have
\begin{equation}
\label{com}R_{_I}=-R_{_I}\andd \LL_\a\sub I;\;\;\a\in R_{_I}.
\end{equation} This  in particular implies that  $I=(I\cap\LL_0)\op\sum_{\a\in
R_{_I}}\LL_\a.$ Now this together with Proposition
\ref{loc5}$(ii)$ implies that $R_{_I}$ is closed. For the last
assertion, one can easily check that $I_c:=\sum_{\a\in
R_{_I}}[\LL_\a,\LL_{-\a}]\op\sum_{\a\in R_{_I}}\LL_\a$ is an ideal
of $I$ and so   we are done as $I$  is simple.
\smallskip

($b$) Since $\LL$ is semisimple, there are an index set $A$ and
simple ideals $\LL^i,$ $i\in A,$  of $\LL$ such that
$\LL=\op_{i\in A}\LL^i.$ For $i\in A,$ use the notation as in the
previous part and set $R_i:=R_{_{\LL^i}},$ since $\LL^i$ is
simple,  part $(a)$ implies  that $\LL^i=\sum_{\a\in
R_i}[\LL_{\a},\LL_{-\a}]\op\sum_{\a\in R_i}\LL_\a.$ Now we are
done as $\LL=\op_{i\in A}\LL^i.$

\smallskip

($c$) Suppose $\a,\b\in R,$ then thanks to   (\ref{com}),   we
have
\begin{equation}\label{closed}\begin{array}{l}
\; [\LL_\a,\LL_\b]\sub\LL_{\a+\b}\cap I\hbox{ if $\a\in  R_{_I}$
or $\b\in R_{_I}$}\\\\ \;[\LL_{\a+\b},\LL_{-\a}]\sub\LL_\b\cap
I\hbox{ if $\a\in R_{_I}$ or $\a+\b\in  R_{_I}$}.\end{array}
\end{equation}
So using  Proposition \ref{loc5}$(ii),$ we have
\begin{equation} \label{clos1}
[\LL_\a,\LL_{\b}]\sub\left\{\begin{array}{ll}
I& \hbox{if $\a,\b\in R_{_I},$}\\
J& \hbox{if $\a,\b\in\D,$}\\
\{0\}& \hbox{if $\a\in R_{_I},\;\b\in\D.$}
\end{array}\right.
\end{equation}
This in particular implies that
\begin{equation}
[\sum_{\a\in R_{_I}}[\LL_\a,\LL_{-\a}]+\sum_{\a\in
R_{_I}}\LL_\a,\sum_{\b\in\D}[\LL_\b,\LL_{-\b}]]=\{0\}.
\label{lc2}\end{equation} Now (\ref{clos1}),(\ref{lc2}) imply that
\begin{eqnarray}
[\LL,J]=[\sum_{\a\in R}\LL_\a+\sum_{\a\in R}[\LL_\a,\LL_{-\a}],J]
&=&[\sum_{\a\in \D}\LL_\a+\sum_{\a\in
\D}[\LL_\a,\LL_{-\a}],J]\nonumber\\
&\sub& [J,J]\sub J\label{ideal}
\end{eqnarray}
which means that $J$ is an ideal of $\LL.$ Also as $\sum_{\a\in
R_{_I}}[\LL_\a,\LL_{-\a}]\sub I,$ one gets that
\begin{equation}\label{dec}\LL=I+J.
\end{equation}

Now suppose that  $x\in I\cap J,$ then $x\in
I\cap\sum_{\a\in\D}[\LL_\a,\LL_{-\a}],$ so by (\ref{lc2}),
\begin{equation}\label{cen1}[x,\sum_{\a\in R_{_I}}\LL_\a\op\sum_{\a\in R_{_I}}[\LL_\a,\LL_{-\a}]]=\{0\}.
\end{equation}
On the other hand since $x\in I\cap\LL_0,$  we have
$$[x,\LL_\a]=\{0\};\;\;\a\in\D.$$
This in particular implies that
$[x,\sum_{\a\in\D}[\LL_\a,\LL_{-\a}]+\sum_{\a\in\D}\LL_\a]=\{0\}.$
Now this together with (\ref{dec}) and  (\ref{cen1})  implies
that $x\in Z(\LL).$  \qed

\begin{cor}
\label{semi-sim} Suppose that $R=R_{int}.$ Set
$\LL_c:=\LL_{0,0}\op\sum_{\a\in R}\LL_\a$ where
$\LL_{0,0}:=\sum_{\a\in R}[\LL_\a,\LL_{-\a}].$ Then
$\LL_{cc}:=\LL_c/Z(\LL_c) $ is a semisimple Lie algebra. In
particular, if $\LL$ is finite dimensional, $\LLc$ is semisimple.
\end{cor}
\pf It follows from  (\ref{con}) that   the canonical projection
map $\LL_c\longrightarrow \LL_{cc}$ restricted to $\sum_{\a\in
R}\LL_\a$ is injective, so we identify $\sum_{\a\in R}\LL_\a$ as a
subspace of $\LL_{cc},$ also it is easy to see that
$(\LL_{cc},\frac{\hh+Z(\LLc)}{Z(\LLc)})$  is an admissible  Lie
algebra whose root system can be identified with  $R.$ More
precisely, we have  the following  weight space decomposition for
$\LL_{cc}$ with respect to $(\hh+Z(\LLc))/{Z(\LLc)}:$
\begin{equation}
\label{llc}\begin{array}{c}\LL_{cc}=(\LL_{cc})_0\op\sum_{\a\in
R}(\LL_{cc})_\a \hbox{ with }\vspace{3mm}\\
 (\LL_{cc})_0=\LL_{0,0}/Z(\LL_c),
(\LL_{cc})_\a=\LL_\a,$ $\a\in R.\end{array}\end{equation}
Moreover we have  $R=R_{int},$ $(\LLcc)_0=\sum_{\a\in
R}[(\LLcc)_\a,(\LLcc)_{-\a}]$ and $Z(\LLcc)=\{0\}.$
 Using
Proposition \ref{semi-sim1}($c$), we get that for any ideal $I$ of
$\LLcc,$ there is an ideal $J$ of $\LLcc$ such that $\LL=I\op J,$
this results in that  $\LLcc$ is semisimple (see \cite[\S
3.5]{J}). Now suppose $\LL$ is finite dimensional and take
$\mathfrak{r}$ to be the solvable radical of $\LLc,$ then
$Z(\LLc)\sub \mathfrak{r}.$ Also as $\mathfrak{r}$ is solvable, we
get  that $\mathfrak{r}/Z(\LLc)$ is a solvable ideal of the
semisimple Lie algebra $\LLcc.$  So $\mathfrak{r}/Z(\LLc)=\{0\},$
i.e., $\mathfrak{r}=Z(\LLc).$ This means that $\LLc$ is a
reductive Lie algebra, now  as $\LLc$ is perfect, one concludes
that $\LLc=[\LLc,\LLc]$ is semisimple. This completes the
proof.\qed

\begin{lem}
\label{rest2} For a  subset  $\D$ of $ R,$ set
$\LL_{_\D}:=\sum_{\a\in\D}[\LL_\a,\LL_{-\a}]\op\sum_{\a\in
\D}\LL_\a$ and $\hh_{_\D}:=\hh\cap\sum_{\a\in
\D}[\LL_\a,\LL_{-\a}].$ Suppose that $\D\sub R_{int}$ is symmetric
and closed and take $\pi_{_\D}:\hh^\star\longrightarrow
\hh_{_\D}^\star$ to be defined by
$\a\mapsto\a_{\mid_{\hh_{_\D}}},$ $\a\in\hh^\star.$ Then
$\pi_{_\D}$ restricted to $\D\cup\{0\}$ is injective. Also
identifying $\D$ with  $\pi_{_\D}(\D),$ we have that
$(\LL_{_\D},\hh_{_\D})$ is an admissible Lie algebra with root
system $\D.$   If moreover   $\D$  is connected with respect to a
fixed  set  of splitting elements, then $\LL_{_\D}$ is  a simple
Lie subalgebra of $\LL.$
\end{lem}
\pf Set $\pi:=\pi_{_\D}.$ Suppose that $\a,\b\in \D\cup\{0\}$ and
$\pi(\a)=\pi(\b).$ We first suppose that $\gamma:=\a-\b\in R,$
then since $\D$ is symmetric and closed, we have  $\gamma\in \D,$
so (\ref{con}) guarantees the existence  of
$t\in\hh\cap[\LL_\gamma,\LL_{-\gamma}]\sub\hh\cap\LL_{_\D}=\hh_{_\D}$
with $\gamma(t)=2,$ so $(\a-\b)(t)=2$ which contradicts the fact
that $\a\mid_{\hh_{_\D}}=\b\mid_{\hh_{_\D}}.$ Next suppose that
$\a-\b\not \in R\cup\{0\},$ then $\a\neq 0$ and $\b\neq0 .$ Since
$\a$ is integrable, there are $e\in\LL_\a,$ $f\in\LL_{-\a}$ such
that $(e,h:=[e,f],f)$ is an $\mathfrak{sl}_2$-triple and $ad_e,$
$ad_f$ act locally nilpotently on $\LL.$ Now since $\a-\b\not\in
R\cup\{0\},$ Proposition \ref{loc5}$(iv)$ implies that
$\b(h)\leq 0.$ This contradicts the fact that $\pi(\a)=\pi(\b)$ as
$h\in\hh\cap[\LL_\a,\LL_{-\a}]\sub\hh_{_\D},$ $\a(h)=2$ and
$\b(h)\leq 0.$ These all together imply that
$\pi\mid_{_{\D\cup\{0\}}}$ is injective. Next we note that since
$\D$ is closed, $\LL_{_\D}$  is a subalgebra of $\LL.$ One also
sees that $\LL_{_\D}$ is an $\hh$-submodule of $\hh$-module $\LL$
admitting the weight space decomposition $\LL_{_\D}=\sum_{\a\in
\D\cup\{0\}}(\LL_{_\D})_\a$ with respect to $\hh$ where
$(\LL_{_\D})_\a=\LL_\a$ for $\a\in \D$ and
$(\LL_{_\D})_0=\sum_{\a\in\D}[\LL_\a,\LL_{-\a}].$ From Lemma
\ref{rest3}, we know that $\LL_{_\D}$ admits a weight space
decomposition
$\LL_{_\D}=\op_{\pi(\a)\in\pi(\D\cup\{0\})}(\LL_{_\D})_{\pi(\a)}$
with respect to $\hh_{_\D}$ where for $\a\in \D\cup\{0\},$
$$(\LL_{_\D})_{\pi(\a)}=\{x\in\LL_{_\D}\mid [h,x]=\a(h)x;\;\forall
h\in\hh_{_\D}\}=\displaystyle{\bigoplus_{\{\b\in
\D\cup\{0\},\pi(\b)=\pi(\a)\}}(\LL_{_\D})_\b}.$$ Now this together
with  the injectivity of $\pi_{\mid_{_{\D\cup\{0\}}}}$  and  the
fact that $(\LL,\hh)$ is an admissible Lie algebra implies that
$(\LL_{_\D},\hh_{_\D})$ is an admissible Lie algebra with root
system $\D.$ The last sentence   follows from Proposition
\ref{semi-sim1}.\qed

\begin{pro}
\label{max-tor} Suppose that $\LL$ is a semisimple Lie algebra and
$\hh$ is a nonzero maximal toral subalgebra of $\LL$ such that
$(\LL,\hh)$ is an admissible Lie algebra with root system
$R=R_{int}.$ Take $I$ to be an index set  such that $\LL=\op_{i\in
I}\LL^i$ where for $i\in I,$ $\LL^i$ is an ideal of $\LL$ which is
simple as a Lie algebra. Set  $R_i:=\{\a\in R\mid \LL_\a\sub
\LL^i\}$ and $\hh_i:=\LL^i\cap\hh,$ $i\in I,$ then

(a) $\hh=\op_{i\in I}\hh_i,$

(b) for $i\in I,$ $\hh_i$ is a maximal toral subalgebra of $\LL^i$
and $(\LL^i,\hh_i)$ is a admissible Lie algebra with root system
$R_i.$
\end{pro}
\pf ($a$) We first mention that for $i,j\in  I,$ we have
\begin{equation}\label{sem}[\LL^i,\LL^j]\sub\d_{i,j}\LL^i.
\end{equation} Using Proposition \ref{semi-sim1}($a$), we get that
for $i\in I,$ $\LL^i=\sum_{\a\in
R_i}[\LL_\a,\LL_{-\a}]\op\sum_{\a\in R_i}\LL_\a,$ and
$\LL_0=\sum_{i\in I}\sum_{\a\in R_i}[\LL_\a,\LL_{-\a}].$ For $i\in
I,$ take $\pi_i:\LL\longrightarrow \LL^i$ to be the natural
projection map and set $T_i$ to be the image of $\hh$ under
$\pi_i.$ Take $T'_i$ to be a subspace of $T_i$ such that
$T_i=(\hh\cap T_i)\op T'_i.$ Let $h\in\hh,$ then  there is a
unique expression $h=\sum_{i\in I}t_i$ with a finitely many
nonzero terms $t_i\in\LL^i,$ $i\in I.$ Now fix $i\in I,$ using
(\ref{sem}), we have for $x\in\LL_\a$ with $\a\in R_i$ that
\begin{equation}\label{ext}
\begin{array}{c}
[\pi_i(h),x]=[t_i,x]=[\sum_{j\in I}t_j,x]=[h,x]=\a(h)x.
\end{array}
\end{equation}
This in particular implies that
\begin{equation*}\label{par2}
\parbox{4in}{if $h,h'\in\hh$  and $\pi_i(h)=\pi_i(h'),$ then   for $\a\in R_i,$ $\a(h)=\a(h')$}
\end{equation*}
which allows us to define $f_{t,\a}$ for $t\in T'_i$ and $\a\in
R_i$ as follows:
\begin{equation}
\label{tor4}
\parbox{4in}{ $f_{t,\a}:=\a(h)$ where $h\in \hh$ is such that $\pi_i(h)=t,$}
\end{equation}
also define $f_{t,\a}:=0$ if $t\in T'_i$ and $\a\in
(R\cup\{0\})\setminus R_i.$
 Next suppose that $\a\in R\cup\{0\}$ and
extend $\a$ to the functional $\a^i\in(\hh\op T'_i)^\star$ defined
as follows:
\begin{equation}
\label{func} \a^i(h+t)=\a(h)+f_{t,\a};\;h\in\hh,\; t\in T'_i.
\end{equation}
Now one can see that $\hh\op T'_i$ is a toral subalgebra of $\LL$
with  corresponding root system $\{\a^i\mid \a\in R\}.$ But
$\hh\op T'_i$ contains the  maximal toral subalgebra $\hh,$
therefore $T'_i=\{0\}$ and so $T_i=T_i\cap\hh$ which implies that
$T_i\sub \hh.$ This gives  that $T_i\sub\hh\cap\LL^i=\hh_i,$ on
the other hand, if $h\in\hh_i,$ then $h=\pi_i(h)\in T_i,$ so we
get that $T_i=\hh_i.$ Now we have $\hh=\op_{i\in I} T_i,$ and so
$\hh=\op_{i\in I}\hh_i.$

$(b)$ Suppose that $i\in I.$ Using Preposition
\ref{semi-sim1}$(a)$ and Lemma \ref{rest2}, we get that
$(\LL^i,\hh_i)$ is an admissible Lie algebra with root system
$R_i.$ Next we note that any toral subalgebra $T_i$ of $\LL_i$
larger than $\hh_i$ would automatically be toral in $\LL$ and
centralizes $\hh_j$ for  $j\in I\setminus\{i\}.$ Now $\hh+T_i$ is
a toral subalgebra of $\LL$ larger than $\hh$ which is a
contradiction. This completes the proof.\qed

\begin{lem}
\label{interesting}  If $(\LL,\hh)$ is a finite dimensional
semi-simple admissible Lie algebra, then $\hh$  is a maximal toral
subalgebra.
\end{lem}
\pf
 Let $\T$ be a maximal toral subalgebra of $\LL$ containing
$\hh.$  We first suppose that $\LL$ is  simple. Since $\T$ is a
toral subalgebra of $\LL,$  $\LL=\sum_{\b\in \T^\star}\LL'_\b$
where for $\b\in \T^\star,$ $\LL'_\b:=\{x\in\LL\mid
[t,x]=\b(t)x;\;\forall t\in \T\}.$ Take $\Phi$ to be the root
system of $\LL$ with respect to $\T.$ Using Example \ref{exa1},
one knows that for $\b\in\Phi,$ there is a unique element
$k_\b\in\T$ satisfying
\begin{equation}\label{property}
\b(k_\b)=2 \hbox{ and $k_\b=[u,v]$ for some
$u\in\LL'_\b,v\in\LL'_{-\b}.$}
\end{equation}

Using Lemma \ref{rest3}, we get  that for $\a\in R\cup\{0\},$
$\LL_\a=\bigoplus_{\b\in A_\a}\LL'_\b$ where
$A_\a:=\{\b\in\T^\star,\;\b\mid_{_\hh}=\a\}.$ Now let $\a\in R$
and $\b\in\Phi$ be such that   $\b\in A_\a.$ Suppose that  $0\neq
x\in\LL'_\b\sub\LL_\a,$ then (\ref{con}) guarantees the existence
of $y\in\LL_{-\a}$ such that $h:=[x,y]\in\hh$ and $(x,h,y)$  is an
$\mathfrak{sl}_2$-triple. Since
$y\in\LL_{-\a}=\bigoplus_{\gamma\in A_\a}\LL'_{-\gamma},$ we have
that $y=\sum_{\gamma\in A_\a}y_{-\gamma}$ with
$y_{-\gamma}\in\LL'_{-\gamma}.$ Now
\begin{eqnarray*}
h=[x,y]=[x,\sum_{\gamma\in A_\a}y_{-\gamma}]= \sum_{\gamma\in
A_\a}[x,y_{-\gamma}]\end{eqnarray*} and so $ h=[x,y_{-\b}] $ as
$h\in\hh\sub\T\sub\LL'_0.$  Also as $h\in \hh,$ $\b(h)=\a(h)=2. $
Now considering (\ref{property}), we get from the uniqueness of
$k_\b$ that $k_\b=h,$ therefore we have proved
\begin{equation}\label{belonh}
k_\b\in\hh;\;\hbox{ for $\b\in\Phi$  with $\b\mid_{_\hh}\neq 0.$}
\end{equation}
Now set $A:=\{\b\in\Phi\mid \b\mid_{_ \hh}\neq 0\}$ and
$I:=\sum_{\b\in A}\LL'_\b\op\sum_{\b\in A}[\LL'_\b,\LL'_{-\b}].$
We claim that $I$ is an ideal of $\LL.$ For this, it is enough to
show that $(\Phi+A)\cap\Phi\sub A.$ Suppose that $\gamma\in\Phi$
and $\b\in A$ are such that
 $\gamma+\b\in\Phi.$ If $\gamma\mid_{_\hh}=0,$ then  $\gamma+\b\mid_{_\hh}\neq0$ and so $\gamma+\b\in A$
 and if
$\gamma\mid_{_\hh}\neq 0,$   then (\ref{lin}) together with
(\ref{belonh}) implies that $k_{\gamma+\b}\in\hh,$ now as
$(\gamma+\b)(k_{\gamma+\b})=2\neq 0,$ we get that
$\gamma+\b\mid_{_\hh}\neq 0,$ i.e., $\gamma+\b\in A.$ This shows
that $I$ is an ideal of $\LL,$ but $\LL$ is simple, so $I=\LL$
which  in particular implies that $A=\Phi.$ Now (\ref{span})
together with (\ref{belonh}) implies that $\T\sub\hh$ and so
$\hh=\T$ is a maximal toral subalgebra. Next suppose that
$\LL=\op_{i=1}^n\LL_i$ where for $1\leq i\leq n,$ $\LL_i$ is a
simple ideal of $\LL.$ So by  Proposition \ref{semi-sim1}($a$) and
Lemma \ref{rest2}, $(\LL_i,\hh\cap\LL_i)$ is a finite dimensional
admissible simple Lie algebra and so using the first part of the
proof, we get that $\hh\cap\LL_i$ is a maximal toral subalgebra of
$\LL_i.$ Now if $\T$ is a maximal toral subalgebra of $\LL$
containing $\hh,$ Example \ref{exa1} gets that $(\LL,\T)$ is a
semi-simple admissible Lie algebra. So Proposition \ref{max-tor}
implies that $\T_i:=\LL_i\cap\T$ is a maximal toral subalgebra of
$\LL_i$ and that $\T=\T_1\op\T_2\op\cdots\op\T_n.$ But
$\hh\cap\LL_i\sub \T_i,$ therefore $\hh\cap\LL_i=\T_i.$ Now  we
have  $\T=\T_1\op\T_2\op\cdots\op\T_n\sub \hh$ and so we are
done.\qed

\begin{pro}\label{simple} (a) Suppose that   $\D$ is a  symmetric closed finite subset of $R_{int}$ such that $\LL_\a$  is finite dimensional for any $\a\in \D,$ then for $\a\in\D,$
there is a unique $h_\a\in\hh$ such that $$\a(h_\a)=2\andd
[x_\a,x_{-\a}]=h_\a\hbox{ for some $x_{\pm\a}\in\LL_{\pm\a}$}.$$
Moreover if  $\a,\b_1,\ldots,\b_m\in\D$ are such that
$\a\in\hbox{span}_\bbbf \{\b_1,\ldots,\b_m\},$ then
$h_\a\in\hbox{span}_\bbbf\{h_{\b_i}\mid1\leq i\leq m\}.$

(b) If  $R=R_{int}$ and all the weight spaces are finite
dimensional, then there is a unique set $\{h_\a\mid \a\in R\}$ of
splitting elements which we refer to as the splitting subset of
$\LL.$
\end{pro}
\pf $(a)$ Set
$\LL_{_\D}:=\sum_{\b\in\D}[\LL_\b,\LL_{-\b}]\op\sum_{\b\in\D}\LL_\b$
and $\hh_{_\D}=\hh\cap\LL_{_\D}.$ Using Propositions \ref{rest2}
and Corollary \ref{semi-sim}, we get that $(\LL_{_\D},\hh_{_\D})$
is a finite dimensional  admissible semi-simple Lie subalgebra
with root system $\D.$ Using  Lemma \ref{interesting}, one  gets
that $\hh_{_\D}$ is a maximal toral subalgebra  of $\LL_{_\D}.$
Now the result follows from  Example \ref{exa1}.

($b$)
 Let $\a\in R$ and take $\D:=R\cap\bbbf\a.$ Using Proposition
\ref{loc5}$(vi)$ together with part ($a$), we are done. \qed

\begin{pro}
\label{bound} Suppose that $(\LL,\hh)$ is a  locally finite
admissible Lie algebra with corresponding root system $R$  such
that all weight spaces are finite dimensional. Take
$\{h_\a\mid\a\in R\}$ to be  the splitting  subset of $\LL$ and
$\v:=\hbox{span}_\bbbf R.$ For $\a\in R,$ define
$\check\a:\v\longrightarrow\bbbf$ by $v\mapsto v(h_\a),$ $v\in\v$
and $s_\a:\v\longrightarrow\v$ by $v\mapsto v-\check\a(v)\a,$
$v\in\v.$ For  $M\sub R,$  set $M_\pm:=M\cup-M$ and take
$\gg_{_M}$ to be the subalgebra of $\LL$ generated by $\cup_{\a\in
M_\pm}\LL_\a,$ finally set $\D_{_M}:=\{\a\in R\mid
\gg_{_M}\cap\LL_\a\neq\{0\}\},$ then we have the followings:

(a) For $M\sub R,$ $\D_{_M}$ is a symmetric closed subset of $R,$
in particular
$\LL_{_{\D_{_M}}}=\sum_{\a\in\D_{_M}}[\LL_\a,\LL_{-\a}]\op\sum_{\a\in\D_{_M}}\LL_\a$
is a Lie subalgebra of $\LL.$

(b) If  $M$ is a  subset of $R$ and $\a,\b\in \D_{_M},$ then
$s_\a(\b)\in\D_{_M}.$ Also if $M$ is a finite subset of $R,$
$\gg_{_M}$ is a finite dimensional Lie subalgebra of $\LL$ and
$\D_{_M}$ is a finite root system in its span where the reflection
based on $\a\in \D_{_M}$ is $s_\a$ restricted to
$\hbox{span}_\bbbf\D_{_M}.$

(c) If $\a,\b\in R,$ then $\b(h_\a)\in\bbbz\cap[-4,4].$

(d) $R$ is a locally finite root system in its span where the
reflection based on $\a\in R$ is $s_\a$ and if $\LL$ is
semi-simple,  the necessary and sufficient condition  for $R$ to
be  irreducible is that $\LL$ is simple. Also for $M\sub R,$
$\D_{_M}$
 is  a closed subsystem of $R$ which is irreducible if $M$ is irreducible.
\end{pro}

\pf $(a)$ We first recall  that as $\LL$ is locally finite, for
any weight vector $x,$ $ad_x$ acts locally nilpotently on $\LL.$
Now we  show that $\D_{_M}$ is symmetric.  One knows that each
element of $\gg_{_M}$ is a sum of elements of the form
$[x_n,\ldots,x_1]$ where $n\in\bbbn\setminus\{0\}$ and for $1\leq
i\leq n,$ $x_i\in\LL_{\a_i},$  for  some $\a_i\in M_{\pm}.$ Now if
$\a\in \D_{_M},$ $\gg_{_M}\cap \LL_\a\neq \{0\},$ so there are
$\a_1,\ldots,\a_n\in M_\pm$ such that $\a=\sum_{i=1}^n\a_i$ and
there are  $x_i\in\LL_{\a_i},$ $1\leq i\leq n,$
 such that $0\neq
x:=[x_n,\ldots,x_1]\in\gg_{_M}\cap \LL_\a.$ We use induction on
$n$ to show that $\D_{_M}$ is symmetric. If $n=1,$ then $\a\in
M_{\pm}$ and so $-\a\in M_{\pm}\sub\D_{_M}.$ Now let
$\a\in\D_{_M}$ be such that there are $n\in\bbbn^{\geq 2},$
$\a_1,\ldots,\a_n\in M_\pm$ and    $x_i\in\LL_{\a_i},$ $1\leq
i\leq n,$  such that $\a=\sum_{i=1}^n\a_i$ and $0\neq
x:=[x_n,\ldots,x_1]\in\gg_{_M}\cap\LL_\a.$ Therefore  $0\neq
[x_{n-1},\ldots,x_1]\in\gg_{_M}\cap \LL_\b$ where
$\b:=\sum_{i=1}^{n-1}\a_i.$ Now  we have
$\b,\a_n,\a=\b+\a_n\in\D_{_M}\sub R.$ Using  the induction
hypothesis, we get that $-\b\in\D_{_M}$ and so there is a nonzero
$y\in\LL_{-\b}\cap\gg_{_M}.$ Now contemplating   Proposition
\ref{loc5}$(ii),$ we have  that $\{0\}\neq
[y,\LL_{-\a_n}]\sub(\gg_{_M}\cap\LL_{-\a_n-\b})$ which shows that
$-\a\in\D_{_M}.$ Next we show that  $\D_{_M}$ is closed. We first
show that
\begin{equation}\label{new1}
\parbox{4.3in}{if $\a\in \D_{_M}$ and $\b\in M_\pm$ are
such that $\a+\b\in R,$ then $\a+\b\in\D_{_M}.$ }\end{equation}
Take $\a\in \D_{_M}$ and $\b\in M_\pm$ to be  such that $\a+\b\in
R$ and x to be a nonzero element of $\gg_{_M}\cap\LL_\a.$ Using
Lemma \ref{loc5}$(ii$), we get that $0\neq
[x,\LL_\b]\sub[\gg_{_M},\gg_{_M}]\sub\gg_{_M}$ and so we are done
as $[x,\LL_\b]\sub\LL_{\a+\b}.$   Now suppose $\a,\b\in \D_{_M}$
and $\a+\b\in R.$  Since $\a,\b\in \D_{_M},$ there is $\{\a_i\mid
1\leq i\leq n\}\cup\{\b_j\mid 1\leq j\leq m\}\sub M_\pm$ such that
$\a=\sum_{i=1}^n\a_i$ and $\b=\sum_{j=1}^m\b_j.$ For $1\leq i\leq
m+n,$ set
$$\gamma_i:=\left\{
\begin{array}{ll}
    \a_i, & \hbox{if } 1\leq i\leq n \\
    \b_{i-n}, & \hbox{if } n+1\leq  i\leq n+m, \\
\end{array}
\right.    $$ then $\gamma_i\in M_\pm,$ $1\leq i\leq m+n$ and
since $\a+\b\in R,$ $\sum_{i=1}^{n+m}\gamma_i\in R.$ Now by
Corollary  \ref{cor0}, there exists a permutation $\pi$ on
$\{1,\ldots,m+n\}$ such that  all partial sums
$\sum_{i=1}^t\gamma_{\pi(i)}\in R\cup\{0\},$ $1\leq t\leq m+n.$
Now using (\ref{new1}) together with an  induction process, we get
that $\a+\b=\sum_{i=1}^{n+m}\gamma_i\in\D_{_M}.$ This means that
$\D_{_M}$ is a closed. This completes the proof.

 $(b)$ Since $\a,\b\in \D_{_M}\sub R,$ by Proposition \ref{loc5}$(i),$
$\check\a(\b)\in\bbbz$ and $s_\a(\b)\in R.$ Now take
$e_\a\in\LL_\a$ and $f_\a\in\LL_{-\a}$ to be  such that
$h_\a=[e_\a,f_\a]$ and $(e_\a,h_\a,f_\a)$ is an
$\mathfrak{sl}_2-$triple, setting  $\theta_\a:=exp(ad_{e_\a
})exp(ad_{-f_\a })exp(ad_{e_\a }),$  we get that $\{0\}\neq
\theta_\a(\LL_\b)\sub(\LL_{s_\a(\b)}\cap\LL_{_{\D_{_M}}})$ and so
$s_\a(\b)\in \D_{_M}.$ Now if $M$ is a finite set, since $\LL$ is
locally finite, $\gg_{_M}$ is a finite dimensional subalgebra of
$\LL$ and so $\D_{_M}$ is finite.  This completes the proof.

$(c)$ Set $M:=\{\a,\b\},$ by part $($b$),$ $\D_{_M}$ is a finite
root system in $\hbox{span}_{\bbbf}\D_{_M}$ where the reflection
based on $\gamma\in\D_{_M}$ is defined by $v\mapsto
v-\check\gamma(v)\gamma,$ $v\in\hbox{span}_{\bbbf}\D_{_M}.$ Using
the finite root system theory, we get  that
 $\{\check\gamma(\eta)\mid
\gamma,\eta\in \D_{_M}\}\sub\bbbz\cap[-4,4],$ in particular
$\b(h_\a)=\check\a(\b)\in\bbbz\cap[-4,4].$

 $(d)$ Using Proposition \ref{loc5}$(i),$ it is  enough to prove that $R$ is locally finite. Suppose
that $M$ is a finite subset of $ R,$ then by part ($a$),
$\D:=\D_{_M}$ is a finite root system in $\hbox{span}_\bbbf\D$
where the reflection based on $\gamma\in\D_{_M}$ is defined by
$v\mapsto v-\check\gamma(v)\gamma$ for $v\in\hbox{span}_\bbbf\D.$
Take $\Pi:=\{\a_1,\ldots,\a_n\}$ to be a base of $\D.$ One knows
that the Cartan matrix  $(\check\a_i(\a_j))$ is an invertible
matrix. So for any choice of $\{k_1,\ldots,k_n\}\sub\bbbf,$ the
following system of equations
$$\sum_{i=1}^n\check\a_j(\a_i)x_i=k_j;\;\; 1\leq j\leq n$$ has a
unique solution. Next suppose that  $(r_1,\ldots,r_n)\in\bbbf^n$
is such that $\eta:=\sum_{i=1}^nr_i\a_i\in R,$ then for any $1\leq
j\leq n,$
$$\sum_{i=1}^nr_i\check\a_j(\a_i)=\sum_{i=1}^nr_i\a_i(h_{_{\a_j}})=\eta(h_{_{\a_j}}).$$
This means that $(r_1,\ldots,r_n)$  is a solution for the
following system of equations
\begin{equation}\label{system}\sum_{i=1}^n\check\a_j(\a_i)x_i=\eta(h_{_{\a_j}});\;\; 1\leq
j\leq n.\end{equation}  But we know that by part ($c$),
$\eta(h_{_{\a_j}})\in [-4,4]\cap\bbbz,$ so there are a finitely
many choice for $\eta(h_{_{\a_j}}),$ $ 1\leq j\leq n.$ This
together with the fact that (\ref{system}) has a unique solution
implies that  there are a finitely many choice for
$(r_1,\ldots,r_n)\in\bbbf^n$ such that $\sum_{i=1}^nr_i\a_i$ is a
root. So $(\hbox{span}_\bbbf M)\cap R=(\hbox{span}_\bbbf\D)\cap
R=(\hbox{span}_\bbbf\{\a_1,\ldots,\a_n\})\cap R$ is a finite set.
Now if $\w$ is a finite  dimensional subspace  of
$\v=\hbox{span}_\bbbf R,$ then there is a finite set $M$ of $R$
such that $\w\sub \hbox{span}_\bbbf M,$ so $\w\cap R\sub
(\hbox{span}_\bbbf M)\cap R$ and so  $\w\cap R$ is a finite set.
This shows that $R$ is a locally finite root system in $\v.$ It
follows easily from  Proposition \ref{semi-sim1}($a$) that if
$\LL$ is semi-simple, $\LL$  is simple if and only if $R$ is
irreducible. Now let $M$ be  a subset of $R,$ parts $(a)$ and
$(b)$ shows that $\D_{_M}$ is a subsystem of $R.$ Next suppose
that $M$ is irreducible but $\D_{_M}$ is not irreducible, then
there is $\a\in\D_{_M}$ such that  $\a$ is connected to none of
the elements of $M,$ in particular,
\begin{equation}\label{conect}
\a(h_\b)=0;\;\;\;\forall\b\in M.
\end{equation}
But   we know there are $\b_1,\ldots,\b_m\in M$ such that
$\a=\sum_{i=1}^m\pm\b_i.$ Take $M':=\{\a,\b_1,\ldots,\b_m\},$ we
know from parts $(a)$ and $(b)$ that  $\D_{_{M'}}$ is a finite
closed  subsystem of $R,$ so Proposition \ref{simple} implies that
$h_\a=\sum_{i=1}^m r_ih_{\b_i}$ for some $r_i\in\bbbf, $ $1\leq
i\leq m.$  This together with (\ref{conect}) implies that
$\a(h_\a)=0$ which is a contradiction. Therefore $\D_{_M}$ is
irreducible. \qed

\bigskip
\begin{pro}
\label{max-tor-ori} Suppose that  $\LL$ is simple and locally
finite such that  all the weight spaces are finite dimensional,
then $\hh$ is a maximal toral subalgebra of $\LL.$
\end{pro}
\pf Take $\{h_\a\mid \a\in R\}$ to be the splitting  subset of
$\LL.$ By Proposition  \ref{bound}$(d)$, $R$ is an irreducible
locally finite root system in its span where the  reflection based
on $\a\in R$ is defined by $v\mapsto v-v(h_\a)\a,$
$v\in\hbox{span}_\bbbf R.$ Now let $T$ be a toral subalgebra of
$\LL$ containing $\hh.$ Take $\Phi$ to be the root system of $\LL$
with respect to $T$ and for $\b\in\Phi\cup\{0\},$ denote by
$\LL'_\b,$ the "weight space" of $\LL$ corresponding to $\b.$ For
$\a\in R\cup\{0\},$ set $A_\a:=\{\b\in\Phi\cup\{0\}\mid
\b\mid_{_\hh}=\a\},$ then by Lemma \ref{rest3}, we have
$\LL_\a=\op_{\b\in A_\a}\LL'_\b.$ Now let $t\in T,$ since $\hh\sub
T$ and $T$ is abelian, $t\in\LL_0,$ but $\LL$ is simple, so
Proposition \ref{semi-sim1}($b$) implies that $t\in\sum_{\a\in
R}[\LL_\a,\LL_{-\a}].$ Therefore there is a finite subset $M$ of
$R$ such that $t\in\sum_{\a\in M}[\LL_\a,\LL_{-\a}].$ Now  take
$\D$ to be  a finite  irreducible closed subsystem of $R$
containing $M$ (see  \cite[Cor. 3.16]{LN}). Using Lemma
\ref{rest2}, we get that
\begin{equation}\label{decom1}\LL_{_\D}=\sum_{\a\in\D}[\LL_\a,\LL_{-\a}]\op\sum_{\a\in \D}\LL_\a=\sum_{\a\in\D}\sum_{\b,\gamma\in
A_\a}[\LL'_\b,\LL'_{-\gamma}]\op\sum_{\a\in \D}\sum_{\b\in
A_\a}\LL'_\b\end{equation} is  an finite
dimensional simple admissible Lie algebra.  Now take
$T_{_\D}:=T\cap\sum_{\a\in\D}[\LL_\a,\LL_{-\a}].$ It is read from
(\ref{decom1}) that $\LL_{_\D}$ has a weight space decomposition
with respect to $T_{_\D},$ in other words $T_{_\D}$ is a toral
subalgebra  of $\LL_{_\D}$ containing $\hh_{_\D}=\hh\cap\LL_{_\D}$
which results in $\hh_{_\D}=T_{_\D}$ contemplating Lemmas \ref{rest2}
and \ref{interesting}. Now we have
$$t\in T\cap \sum_{\a\in M}[\LL_\a,\LL_{-\a}]\sub T\cap
\sum_{\a\in\D }[\LL_\a,\LL_{-\a}]=T_{_\D}=\hh_{_\D}\sub\hh.$$ This
completes the proof.\qed

\medskip

Now one can use Propositions \ref{max-tor-ori}, \ref{bound}($d$)
and  Lemma \ref{rest2}  to get the following theorem:
\begin{thm}\label{direct}
Suppose that $(\LL,\hh)$ is a locally finite  simple admissible
Lie algebra whose weight spaces are finite dimensional. Take $R$
to be the root system of $\LL$ with respect to $\hh. $ Then $R$ is an irreducible  locally finite root system. Next take
$\{R_i\mid i\in I\}$  to be the class of finite irreducible closed
subsystems of $R.$  For $i\in I,$ set $\LL_i:=\sum_{\a\in
R_i}[\LL_\a,\LL_{-\a}]\op\sum_{\a\in R_i}\LL_\a$ and
$\hh_i:=\hh\cap\LL_i.$ Then  $\hh$ is a maximal toral subalgebra of $\LL$ and for $i\in I,$  $\hh_i$ is a maximal
toral subalgebra of $\LL_i,$ also  $(\LL_i,\hh_i)$ is  a finite
dimensional simple admissible Lie algebra with root system $R_i.$
Moreover $\{R_i\mid i\in I\}$ and  $\{\LL_i\mid i\in I\}$ are
directed systems with respect to inclusion,  $R$ is a direct limit
of $\{R_i\mid i\in I\}$ and $\LL$ is a direct limit  of
$\{\LL_i\mid i\in I\}.$
\end{thm}


\begin{thebibliography}{99}
\bibitem[ABG]{ABG} B.N. Allison, G. Benkart and Y. Gao,
{\it Lie algebras graded by the root systems $BC_r,$ $r\geq2,$}
 Mem. Amer. Math. Soc.  {\bf 158 } (2002),  no. 751, x+158.


\bibitem[BZ]{BZ}  G. Benkart and E. Zelmanov, {\it Lie algebras graded by finite root systems and intersection matrix
algebras},
   Invent. Math.  {\bf 126}  (1996),  no. 1, 1--45.


\bibitem[BM]{BM} S. Berman and R. Moody, {\it Lie algebras graded by finite root systems and the intersection matrix algebras of
 Slodowy,} Invent. Math.  {\bf 108}  (1992),  no. 2, 323--347.




\bibitem[H]{H} J.E. Humphreys, {\it Introduction to Lie algebras and
representation theory}, Spring Verlag, New York, 1972.
\bibitem[J]{J} N. Jacobson, {\it Basic algebra II}, Second ed., W. H. Freeman and Company, 1989.

\bibitem[K]{K} V. Kac, {\it Infinite dimensional Lie algebras,} third ed., Cambridge University Press, 1990.

\bibitem[LN]{LN} O. Loos and  E. Neher, {\it Locally finite root systems}, Mem. Amer. Math. Soc. {\bf 171} (2004), no. 811, x+214.



\bibitem[MP]{MP} R. Moody and A. Pianzola, {\it Lie algebras with
triangular decomposition,} A Wiley-Interscience Publication, New
York, 1995.

\bibitem[NS]{NS} K.H. Neeb and N. Stumme, {\it  The classification of locally finite split simple Lie algebras}, J. Reine angew. Math. {\bf 533} (2001), 25-53.


\bibitem[N]{N} E. Neher,  {\it Lectures on root-graded and extended affine Lie algebras}, preprint.

\bibitem[Se]{Se} G.B. Seligman, {\it Rational methods in Lie algebras, } M. Dekker Lect. Notes in
pure and appl. math. {\bf 17}, New York, 1976.

\bibitem[St]{St} N. Stumme, {\it  The Structure  of locally finite split simple Lie algebras}, J. Algebra {\bf 220} (1990). no.2, 664--693.


\end{thebibliography}
\end{document}